\documentclass[12pt]{article}
\usepackage{amssymb}
\usepackage{latexsym}
\usepackage{amscd}

\newtheorem{definition}{Definition}[section]
\newtheorem{theorem}[definition]{Theorem}
\newtheorem{lemma}[definition]{Lemma}
\newtheorem{corollary}[definition]{Corollary}
\newtheorem{remark}[definition]{Remark}
\newtheorem{example}[definition]{Example}
\newtheorem{conjecture}[definition]{Conjecture}
\newtheorem{problem}[definition]{Problem}
\newtheorem{note}[definition]{Note}
\newtheorem{assumption}[definition]{Assumption}
\newtheorem{proposition}[definition]{Proposition}

\def\R{\mathbb R}
\def\C{\mathbb C}

\def\Z{\mathbb Z}

\begin{document}

\title{\bf Distance-regular graphs and \\
the $q$-tetrahedron algebra}

\author{
Tatsuro Ito{\footnote{
Department of Computational Science,
Faculty of Science,
Kanazawa University,
Kakuma-machi,
Kanazawa 920-1192, Japan
}}
{\footnote{Supported in part by JSPS grant 18340022}} 
$\;$ and
Paul Terwilliger{\footnote{
Department of Mathematics, University of
Wisconsin, 480 Lincoln Drive, Madison WI 53706-1388 USA}
}}
\date{}

\maketitle

\centerline{\large In honor of Eiichi Bannai on his 60th Birthday}

\begin{abstract}
Let $\Gamma$ denote a distance-regular graph with
classical parameters $(D,b,\alpha,\beta)$ and
$b\not=1$, $\alpha=b-1$. The condition on
$\alpha$ implies that $\Gamma$ is formally self-dual.
For $b=q^2$ we use the adjacency matrix and dual adjacency matrix
to obtain an action of the $q$-tetrahedron algebra
$\boxtimes_q$
on the standard module of $\Gamma$.
We describe four algebra homomorphisms into 
$\boxtimes_q$
from
the quantum affine algebra 
$U_q({\widehat{\mathfrak{sl}}_2})$; 
using these 
we pull back the above
$\boxtimes_q$-action
to obtain four actions of 
$U_q({\widehat{\mathfrak{sl}}_2})$ 
on the standard module of $\Gamma$.

 \end{abstract}

\section{Introduction}
In \cite{HT}
B. Hartwig and the second author
gave a presentation of the three-point $\mathfrak{sl}_2$ loop
algebra via generators and relations. To obtain this
presentation they defined a Lie algebra $\boxtimes$ by
generators and relations, and displayed an isomorphism
 from $\boxtimes$ to
the three-point $\mathfrak{sl}_2$ loop algebra.
The algebra $\boxtimes$ has essentially six generators, and it is natural
to identify these with the six edges of a tetrahedron.
For each face of the
tetrahedron the three 
surrounding edges
form a basis for a subalgebra of $\boxtimes$
that is isomorphic to $\mathfrak{sl}_2$ \cite[Corollary 12.4]{HT}. 
Any five of the six edges of the tetrahedron generate a subalgebra of
$\boxtimes$ 
that is isomorphic to the $\mathfrak{sl}_2$ loop algebra
 \cite[Corollary 12.6]{HT}.
Each pair of opposite
edges of the tetrahedron generate a subalgebra of
$\boxtimes$ that is isomorphic
to the Onsager algebra \cite[Corollary 12.5]{HT}.
Let us call these Onsager subalgebras. Then 
$\boxtimes
$
is the direct sum of its three Onsager subalgebras
 \cite[Theorem 11.6]{HT}.
In \cite{E} Elduque found an 
attractive decomposition of 
$\boxtimes$ into a direct sum of three 
abelian subalgebras, and he
showed how these subalgebras 
are related to the Onsager subalgebras. 
In \cite{PT} Pascasio and the second author give an
action of $\boxtimes$ on the standard module of a Hamming graph.
 In \cite{Br} Bremner
obtained 
the universal central extension of 
the three-point $\mathfrak{sl}_2$ loop algebra.
By modifying the defining relations for $\boxtimes$, 
Benkart and the second author obtained  
a presentation for this
extension by generators and relations
 \cite{BT}. 
In \cite{Ha} Hartwig obtained the 
irreducible finite-dimensional $\boxtimes$-modules
over an algebraically closed field with characteristic 0.

\medskip
\noindent 
In \cite{qtet} we introduced a quantum analog of
$\boxtimes$ which we call 
$\boxtimes_q$.
We defined $\boxtimes_q$ using generators and relations.
We showed how $\boxtimes_q$ is related to the
quantum group
$U_q(\mathfrak{sl}_2)$ in roughly the same
way that $\boxtimes $ is related to $\mathfrak{sl}_2$ \cite[Proposition 7.4]{qtet}.
We showed how $\boxtimes_q$ is related to
the 
$U_q(\mathfrak{sl}_2)$ loop algebra
in roughly the same way that
$\boxtimes$ is related to the $\mathfrak{sl}_2$ loop algebra 
\cite[Proposition 8.3]{qtet}.
In \cite{NN} we considered an algebra
${\mathcal A}_q$ on two generators
subject to the cubic $q$-Serre relations.
${\mathcal A}_q$ is often called the
{\it positive part  of
$U_q({\widehat{\mathfrak{sl}}_2})$}.
We showed how $\boxtimes_q$ is related to
${\mathcal A}_q$ in roughly the same way that
$\boxtimes$ is related to the Onsager algebra \cite[Proposition 9.4]{qtet}.
In \cite{qtet} and \cite{qinv}
we described the finite-dimensional irreducible
$\boxtimes_q$-modules under the assumption
 that $q$ is not a root of 1, and the
 underlying field is algebraically closed.

\medskip
\noindent
 In the
present paper we
consider a distance-regular graph $\Gamma$ that has classical parameters
$(D,b,\alpha,\beta)$ and $b\not=1$, $\alpha=b-1$. The condition on
$\alpha$ implies that $\Gamma$ is formally self-dual \cite[p.~71]{bcn}.
For $b=q^2$ we use the adjacency matrix and dual adjacency matrix to
construct an action of $\boxtimes_q$ on the standard module of $\Gamma$.
We describe four algebra homomorphisms from
$U_q({\widehat{\mathfrak{sl}}_2})$ 
to $\boxtimes_q$; using these homomorphisms
we pull back the above
$\boxtimes_q$-action
to obtain four actions of 
$U_q({\widehat{\mathfrak{sl}}_2})$ 
on the standard module of $\Gamma$.
Several well-known families of distance-regular graphs satisfy the
above parameter restriction; for instance the
bilinear forms graph \cite[p.~280]{bcn},
the alternating forms graph \cite[p.~282]{bcn},
the Hermitean forms graph \cite[p.~285]{bcn},
the quadratic forms graph \cite[p.~290]{bcn},
the affine $E_6$ graph \cite[p.~340]{bcn},
and the extended ternary Golay code graph \cite[p.~359]{bcn}.

\medskip
\noindent All of the original results in this paper are about
distance-regular graphs. However, in order to motivate things
and develop some machinery, we will initially discuss 
$\boxtimes_q$ and its relationship to certain quantum groups.
The paper is organized as follows. In Section 2 we
define $\boxtimes_q$ and mention a few of its properties.
In Section 3 we recall how $\boxtimes_q$ is related to
$U_q(\mathfrak{sl}_2)$. 
In Section 4 we discuss how $\boxtimes_q$ is related to
$U_q({\widehat{\mathfrak{sl}}_2})$.
In Section 5 we recall how $\boxtimes_q$ is related to
${\mathcal A}_q$. In Section 6 we discuss the
finite-dimensional irreducible $\boxtimes_q$-modules.
In Section 7 we consider a distance-regular graph
$\Gamma $ and discuss its basic properties.
In Section 8 we impose a parameter restriction on
$\Gamma$ needed to construct our
$\boxtimes_q$-module.
In Sections 9, 10 we define some matrices that will be used
 to construct
our $\boxtimes_q$-module.
In Section 11 we display an action of $\boxtimes_q$ on
the standard module of $\Gamma$; Theorem 11.1 is the
main result of the paper.
In Section 12 we discuss how the above $\boxtimes_q$-action is
related to the subconstituent algebra of $\Gamma$. 
In Section 13 we give some suggestions for further research.

\medskip
\noindent Throughout the paper $\C$ denotes the field of
complex numbers.

\section{The $q$-tetrahedron algebra $\boxtimes_q$}

\noindent In this section we recall the 
$q$-tetrahedron algebra. 
We fix a nonzero scalar $q \in \C$ such that $q^2\not=1$ and define
\begin{eqnarray*}
\lbrack n \rbrack_q = \frac{q^n-q^{-n}}{q-q^{-1}},
\qquad \qquad n = 0,1,2,\ldots 
\label{eq:nbrack}
\end{eqnarray*}
We let $\Z_4 = \Z/4\Z$ denote the cyclic group of order 4.

\begin{definition} \rm \cite[Definition 10.1]{qtet}
\label{def:qtet}
Let $\boxtimes_q$ denote the unital associative $\C$-algebra that has
generators 
\begin{eqnarray*}
\lbrace x_{ij}\;|\; i,j \in \Z_4,\;j-i=1 \;\mbox{or} \;j-i=2\rbrace
\end{eqnarray*}
and the following relations:
\begin{enumerate}
\item For $i,j\in \Z_4$ such that $j-i=2$,
\begin{eqnarray*}
x_{ij}x_{ji} = 1.
\label{eq:qrel0}
\end{eqnarray*}
\item For $h,i,j\in \Z_4$ such that the pair $(i-h,j-i)$ is one of
$(1,1), (1,2), (2,1)$,
\begin{eqnarray*}
\frac{qx_{hi}x_{ij}-q^{-1}x_{ij}x_{hi}}{q-q^{-1}}=1.
\label{eq:qrel1}
\end{eqnarray*}
\item For $h,i,j,k\in \Z_4$ such that $i-h=j-i=k-j=1$,
\begin{eqnarray}
\label{eq:qserre}
x_{hi}^3x_{jk} -
\lbrack 3 \rbrack_q
x_{hi}^2x_{jk}x_{hi} +
\lbrack 3 \rbrack_q
x_{hi}x_{jk}x_{hi}^2- 
x_{jk}x_{hi}^3=0. 
\end{eqnarray}
\end{enumerate}
We call $\boxtimes_q$ the
{\it $q$-tetrahedron algebra} or ``$q$-tet'' for short.
\end{definition}
\begin{note}\rm
The equations (\ref{eq:qserre}) are the cubic $q$-Serre relations
\cite[p.~10]{lusztig}.
\end{note}
\noindent
We make some observations.

\begin{lemma}
\label{lem:rho}
{\rm \cite[Lemma 6.3]{qtet} }
There exists a $\C$-algebra automorphism $\varrho$ of $\boxtimes_q$
that sends each generator $x_{ij}$ to $x_{i+1,j+1}$.
Moreover 
 $\varrho^4=1$.
\end{lemma}


\begin{lemma}
\label{lem:flip}
{\rm \cite[Lemma 6.5]{qtet}}
There exists a $\C$-algebra automorphism of 
 $\boxtimes_q$
that sends each generator $x_{ij}$ to $-x_{ij}$.
\end{lemma}

\section{The algebra $U_q(\mathfrak{sl}_2)$}

\noindent In this section we recall how the algebra
$\boxtimes_q$ is related to 
$U_q(\mathfrak{sl}_2)$.
We start with a definition.

\begin{definition} 
\label{def:uq}
\rm
\cite[p.~122]{Kassel}
Let $U_q(\mathfrak{sl}_2)$
denote the unital associative $\C$-algebra 
with
generators $K^{\pm 1}$, $e^{\pm}$
and the following relations:
\begin{eqnarray*}
KK^{-1} &=& 
K^{-1}K =  1,
\label{eq:buq1}
\\
Ke^{\pm}K^{-1} &=& q^{\pm 2}e^{\pm},
\label{eq:buq2}
\\
\lbrack e^+,e^-\rbrack  &=& \frac{K-K^{-1}}{q-q^{-1}}.
\label{eq:buq4}
\end{eqnarray*}
\end{definition}


\noindent The following presentation of
$U_q(\mathfrak{sl}_2)$ will be useful.

\begin{lemma}
\label{thm:uq2}
{\rm \cite[Theorem 2.1]{equit1}}
The algebra
$U_q(\mathfrak{sl}_2)$ 
 is isomorphic to
the unital associative $\C$-algebra 
with
generators 
$x^{\pm 1}$, $y$, $z$
and the following relations:
\begin{eqnarray*}
xx^{-1} = 
x^{-1}x &=&  1,
\label{eq:2buq1}
\\
\frac{qxy-q^{-1}yx}{q-q^{-1}}&=&1,
\label{eq:2buq2}
\\
\frac{qyz-q^{-1}zy}{q-q^{-1}}&=&1,
\label{eq:2buq3}
\\
\frac{qzx-q^{-1}xz}{q-q^{-1}}&=&1.
\label{eq:2buq4}
\end{eqnarray*}
An isomorphism with the presentation in Definition
\ref{def:uq} is given by:
\begin{eqnarray*}
\label{eq:iso1}
x^{{\pm}1} &\mapsto & K^{{\pm}1},\\
\label{eq:iso2}
y &\mapsto & K^{-1}+e^-, \\
\label{eq:iso3}
z &\mapsto & K^{-1}-K^{-1}e^+q(q-q^{-1})^2.
\end{eqnarray*}
The inverse of this isomorphism is given by:
\begin{eqnarray*}
\label{eq:iso1inv}
K^{{\pm}1} &\mapsto & x^{{\pm}1},\\
\label{eq:iso2inv}
e^- &\mapsto & y-x^{-1}, \\
\label{eq:iso3inv}
e^+ &\mapsto & (1-xz)q^{-1}(q-q^{-1})^{-2}.
\end{eqnarray*}
\end{lemma}


\begin{proposition}
\label{lem:uqinj}
{\rm \cite[Proposition 7.4]{qtet}}
For $i \in \Z_4$ there exists a $\C$-algebra homomorphism
from
$U_q(\mathfrak{sl}_2)$ to
 $\boxtimes_q$ that sends 
\begin{eqnarray*}
x\mapsto x_{i,i+2},
\quad 
x^{-1}\mapsto x_{i+2,i},\quad
y\mapsto x_{i+2,i+3},
\quad 
z\mapsto x_{i+3,i}.
\end{eqnarray*}
\end{proposition}

\section{The quantum affine algebra
$U_q(\widehat{ \mathfrak{sl}}_2)$
}

\noindent In this section we consider how
$\boxtimes_q$ is related to the
quantum affine algebra
$U_q(\widehat{ \mathfrak{sl}}_2)$.
We start with a definition.

\begin{definition} 
\label{def:qauq}
\rm
\cite[p.~266]{charp} 
The quantum affine algebra
$U_q(\widehat{ \mathfrak{sl}}_2)$
is the unital associative $\C$-algebra 
with
generators $K^{\pm 1}_i$,
 $e^{\pm}_i$,
$i\in \lbrace 0,1\rbrace $
and the following relations:
\begin{eqnarray*}
K_iK^{-1}_i &=& K^{-1}_iK_i=1,
\label{eq:qauq1}
\\
K_0K_1&=& K_1K_0,
\label{eq:qauq2}
\\
K_ie^{\pm}_iK^{-1}_i &=& q^{{\pm}2}e^{\pm}_i,
\label{eq:qauq3}
\\
K_ie^{\pm}_jK^{-1}_i &=& q^{{\mp}2}e^{\pm}_j, \qquad i\not=j,
\label{eq:qauq4}
\\
\lbrack e^+_i, e^-_i\rbrack &=& {{K_i-K^{-1}_i}\over {q-q^{-1}}},
\label{eq:qauq5}
\\
\lbrack e^{\pm}_0, e^{\mp}_1\rbrack &=& 0,
\label{eq:qauq6}
\end{eqnarray*}
\begin{eqnarray*}
(e^{\pm}_i)^3e^{\pm}_j -  
\lbrack 3 \rbrack_q (e^{\pm}_i)^2e^{\pm}_j e^{\pm}_i 
+\lbrack 3 \rbrack_q e^{\pm}_ie^{\pm}_j (e^{\pm}_i)^2 - 
e^{\pm}_j (e^{\pm}_i)^3 =0, \qquad i\not=j.
\label{eq:qauq7}
\end{eqnarray*}
\end{definition}

\noindent The following presentation of
$U_q(\widehat{ \mathfrak{sl}}_2)$
 will be useful.

\begin{theorem}
\label{thm:qa2} 
{\rm (\cite[Theorem 2.1]{tdanduq},
\cite{equit2})}
The quantum affine algebra
$U_q(\widehat{ \mathfrak{sl}}_2)$
 is isomorphic to
the unital associative $\C$-algebra 
with
generators $x_i^{\pm 1}$, $y_i$, $z_i$, $i\in \lbrace 0,1\rbrace $
and the following relations:
\begin{eqnarray*}
x_ix^{-1}_i = x^{-1}_ix_i &=&1,\\
x_0x_1 \;\;\mbox{is central},
\label{eq:qabuq2}
\\
\frac{q x_iy_i-q^{-1}y_ix_i}{q-q^{-1}} &=& 1,
\label{eq:qabuq3}
\\
\frac{q y_iz_i-q^{-1}z_iy_i}{q-q^{-1}} &=& 1,
\label{eq:qabuq4}
\\
\frac{q z_ix_i-q^{-1}x_iz_i}{q-q^{-1}} &=& 1,
\label{eq:qabuq5}
\\
\frac{q z_iy_j-q^{-1}y_jz_i}{q-q^{-1}} &=&x^{-1}_0x^{-1}_1,
\qquad i\not=j,
\label{eq:qabuq6}
\end{eqnarray*}
\begin{eqnarray*}
y_i^3y_j -  
\lbrack 3 \rbrack_q y_i^2y_j y_i 
+\lbrack 3 \rbrack_q y_iy_j y_i^2 - 
y_j y_i^3 =0, \qquad i\not=j,
\label{eq:qabuq7}
\\
z_i^3z_j -  
\lbrack 3 \rbrack_q z_i^2z_j z_i 
+\lbrack 3 \rbrack_q z_iz_j z_i^2 - 
z_j z_i^3 =0, \qquad i\not=j.
\label{eq:qabuq8*}
\end{eqnarray*}
An isomorphism with the presentation in Definition
\ref{def:qauq} is given by:
\begin{eqnarray*}
\label{eq:qaiso1}
x^{\pm 1}_i &\mapsto & K^{\pm 1}_i,\\
\label{eq:qaiso2}
y_i &\mapsto & K^{-1}_i+e^{-}_i, \\
\label{eq:qaiso3}
z_i &\mapsto & K^{-1}_i-K^{-1}_ie^{+}_iq(q-q^{-1})^2.
\end{eqnarray*}
The inverse of this isomorphism is given by:
\begin{eqnarray*}
\label{eq:qliso1inv}
K^{\pm 1}_i &\mapsto & x^{\pm 1}_i,\\
\label{eq:qliso2inv}
e^-_i &\mapsto & y_i-x^{-1}_i, \\
\label{eq:qliso3inv}
e^+_i &\mapsto & (1-x_i z_i)q^{-1}(q-q^{-1})^{-2}.
\end{eqnarray*}
\end{theorem}

\begin{proposition}
\label{prop:uqhom}
For $i \in \Z_4$ there exists a $\C$-algebra homomorphism
from 
$U_q(\widehat{ \mathfrak{sl}}_2)$
to $\boxtimes_q$
that sends
\begin{eqnarray*}
&&
x_1 \mapsto x_{i,i+2},\quad
x^{-1}_1 \mapsto x_{i+2,i},\quad
y_1\mapsto x_{i+2,i+3},\quad
z_1 \mapsto x_{i+3,i},
\\
&&x_0 \mapsto x_{i+2,i}, \quad 
x^{-1}_0 \mapsto x_{i,i+2}, \quad
y_0 \mapsto x_{i,i+1},\quad
z_0 \mapsto x_{i+1,i+2}.
\end{eqnarray*}
\end{proposition}
\noindent {\it Proof:} 
Compare the defining relations for
$U_q(\widehat{ \mathfrak{sl}}_2)$
 given in
Theorem
\ref{thm:qa2} with the relations
in Definition
\ref{def:qtet}.
\hfill $\Box $

\section{The algebra ${\mathcal A}_q$}

In this section we recall how $\boxtimes_q$ is 
related to the algebra ${\mathcal A}_q$.
We start with a definition.

\begin{definition}
\label{defa}
\rm
Let ${\mathcal A}_q$ denote the unital associative $\C$-algebra
defined by generators $x,y$
and relations
\begin{eqnarray*}
\label{eq:fqs1}
x^3y-\lbrack 3\rbrack_q x^2yx
+\lbrack 3\rbrack_q xyx^2
-yx^3&=&0,
\\
\label{eq:fqs2}
y^3x-\lbrack 3\rbrack_q y^2xy
+\lbrack 3\rbrack_q yxy^2
-xy^3&=&0.
\end{eqnarray*}
\end{definition}

\begin{definition} \rm
Referring to Definition
\ref{defa}, we call $x,y$ the {\it standard generators}
for 
${\mathcal A}_q$.
\end{definition}

\begin{note}\rm
{\rm \cite[Corollary 3.2.6]{lusztig}}
The algebra ${\mathcal A}_q$ is often called the
{\it positive part of 
$U_q(\widehat{ \mathfrak{sl}}_2)$.}
\end{note}

\begin{proposition}
\label{lem:aqinj}
{\rm \cite[Proposition 9.4]{qtet}}
For $i \in \Z_4$ there exists a homomorphism of $\C$-algebras
from ${\mathcal A}_q$ to $\boxtimes_q$
that sends the standard generators
$x, y$ to $x_{i,i+1}, x_{i+2,i+3}$ respectively.
\end{proposition}

\section{The finite-dimensional
irreducible $\boxtimes_q$-modules}

\medskip
\noindent In this section we recall how  
 the finite-dimensional irreducible modules
for 
$\boxtimes_q$ and ${\mathcal A}_q$ are
related.
We start with some comments.
Let $V$ denote a 
finite-dimensional vector space over $\C$. 
A linear transformation $A:V\to V$ is said
to be {\it nilpotent} whenever
there exists a positive integer $n$ such that
$A^n=0$. 
Let $V$ denote a finite-dimensional irreducible 
${\mathcal A}_q$-module.
This  module is called {\it NonNil}
whenever the standard generators $x, y$ are not nilpotent
on $V$ \cite[Definition 1.3]{NN}.
Assume $V$ is
NonNil.
Then
by \cite[Corollary 2.8]{NN} 
 the standard generators $x,y$ are
semisimple on
$V$. Moreover
there exist an integer $d\geq 0$ and nonzero scalars $\alpha, \alpha^* \in \C$
such that the set of distinct eigenvalues of $x$ (resp. $y$) on $V$
is
$\lbrace \alpha q^d,  
\alpha q^{d-2},
\ldots, \alpha q^{-d}
\rbrace$
(resp. 
$\lbrace \alpha^* q^d,  
\alpha^* q^{d-2},
\ldots, \alpha^* q^{-d}
\rbrace$).
We call the ordered pair $(\alpha,\alpha^*)$ the {\it type}
of $V$.
Replacing $x,y$ by $x/\alpha, y/\alpha^*$ the type becomes $(1,1)$.
Now let $V$ denote a finite-dimensional irreducible 
$\boxtimes_q$-module. By \cite[Theorem 12.3]{qtet} each 
generator $x_{ij}$ is semisimple on $V$. Moreover
there exist an
integer $d\geq 0$ and a scalar $\varepsilon \in \lbrace 1,-1\rbrace$
such that for each generator $x_{ij}$ the
set of distinct eigenvalues on  
$V$ is 
$\lbrace \varepsilon q^d, 
 \varepsilon q^{d-2}, \ldots,
 \varepsilon q^{-d}\rbrace$.  
We call $\varepsilon $ the {\it type} of $V$.
Replacing each generator $x_{ij}$ by $\varepsilon x_{ij}$ the
type becomes 1. 
The finite-dimensional
irreducible modules for 
$\boxtimes_q$ and ${\mathcal A}_q$ are related according
to the following two theorems and subsequent remark.

\begin{theorem}
\label{thm:2}
{\rm \cite[Theorem 10.3]{qtet}}
Let $V$ denote a 
finite-dimensional
irreducible 
$\boxtimes_q$-module of type 1.
Then there exists
a unique
${\mathcal A}_q$-module 
structure on 
$V$ such that the standard generators
$x$ and  $y$ act as $x_{01}$  
and $x_{23}$ respectively.
This 
${\mathcal A}_q$-module 
is irreducible, NonNil, and 
 type  $(1,1)$.
\end{theorem}

\begin{theorem}
\label{thm:1}
{\rm \cite[Theorem 10.4]{qtet}}
Let $V$ denote a
NonNil 
finite-dimensional
irreducible
${\mathcal A}_q$-module of type $(1,1)$.
Then there exists
a unique
$\boxtimes_q$-module
structure on 
$V$ such that the standard generators
$x$ and  $y$ act as $x_{01}$  
and $x_{23}$ respectively.
This 
$\boxtimes_q$-module structure
is irreducible and type $1$.
\end{theorem}

\begin{remark}
\rm \cite[Remark 10.5]{qtet}
Combining Theorem
\ref{thm:2} and
Theorem
\ref{thm:1}
we obtain a bijection between the following two sets:
\begin{enumerate}
\item the isomorphism classes of 
finite-dimensional 
irreducible 
$\boxtimes_q$-modules
of type $1$;
\item the isomorphism classes of NonNil
finite-dimensional 
irreducible 
${\mathcal A}_q$-modules of type $(1,1)$.
\end{enumerate}
\end{remark}

\section{Distance-regular graphs; preliminaries}
We now turn our attention to distance-regular graphs.
After a brief review of their basic properties we 
consider a special type said to be formally self-dual with
classical parameters. From such a distance-regular graph
we will obtain a $\boxtimes_q$-module.

\medskip
\noindent
We now review some definitions and basic concepts concerning distance-regular
graphs.
For more information we refer the reader to 
\cite{bannai,bcn,godsil,terwSub1}.

\medskip
\noindent
Let $X$ denote a nonempty  finite  set.
Let $\hbox{Mat}_X(\C)$ 
denote the $\C$-algebra
consisting of all matrices whose rows and columns are indexed by $X$
and whose entries are in $\C  $. Let
$V=\C^X$ denote the vector space over $\C$
consisting of column vectors whose 
coordinates are indexed by $X$ and whose entries are
in $\C$.
We observe
$\hbox{Mat}_X(\C)$ 
acts on $V$ by left multiplication.
We call $V$ the {\it standard module}.
We endow $V$ with the Hermitean inner product $\langle \, , \, \rangle$ 
that satisfies
$\langle u,v \rangle = u^t\overline{v}$ for 
$u,v \in V$,
where $t$ denotes transpose and $\overline{\phantom{v}}$
denotes complex conjugation.
For all $y \in X,$ let $\hat{y}$ denote the element
of $V$ with a 1 in the $y$ coordinate and 0 in all other coordinates.
We observe $\{\hat{y}\;|\;y \in X\}$ is an orthonormal basis for $V.$

\medskip
\noindent
Let $\Gamma = (X,R)$ denote a finite, undirected, connected graph,
without loops or multiple edges, with vertex set $X$ and 
edge set
$R$.   
Let $\partial $ denote the
path-length distance function for $\Gamma $,  and set
$D := \mbox{max}\{\partial(x,y) \;|\; x,y \in X\}$.  
We call $D$  the {\it diameter} of $\Gamma $.
For an integer $k\geq 0$ we say that $\Gamma$ is {\it regular with
valency $k$} whenever each vertex of $\Gamma$ is adjacent to
exactly $k$ distinct vertices of $\Gamma$.
 We say that $\Gamma$ is {\it distance-regular}
whenever for all integers $h,i,j\;(0 \le h,i,j \le D)$ 
and for all
vertices $x,y \in X$ with $\partial(x,y)=h,$ the number
\begin{eqnarray*}
p_{ij}^h = |\{z \in X \; |\; \partial(x,z)=i, \partial(z,y)=j \}|
\end{eqnarray*}
is independent of $x$ and $y.$ The $p_{ij}^h$ are called
the {\it intersection numbers} of $\Gamma.$ 
We abbreviate $c_i=p^i_{1,i-1}$ $(1 \leq i \leq D)$,
$b_i=p^i_{1,i+1}$ $(0 \leq i \leq D-1)$,
$a_i=p^i_{1i}$ $(0 \leq i \leq D)$.

\medskip
\noindent
For the rest of this paper we assume  $\Gamma$  
is  distance-regular; to avoid trivialities we always
assume  $D\geq 3$.
Note that $\Gamma$ is regular with valency $k=b_0$. Moreover
 $k=c_i+a_i+b_i$ for $0 \leq i \leq D$, where $c_0=0$ and
$b_D=0$. 

\medskip
\noindent 
We mention a fact for later use.
By the triangle inequality, for $0 \leq h,i,j\leq D$ we have
$p^h_{ij}= 0$
(resp. 
$p^h_{ij}\not= 0$) whenever one of $h,i,j$ is greater than
(resp. equal to) the sum of the other two.

\medskip
\noindent 
We recall the Bose-Mesner algebra of $\Gamma.$ 
For 
$0 \le i \le D$ let $A_i$ denote the matrix in $\hbox{Mat}_X(\C)$ with
$(x,y)$-entry
$$
{(A_i)_{xy} = \cases{1, & if $\partial(x,y)=i$\cr
0, & if $\partial(x,y) \ne i$\cr}} \qquad (x,y \in X).
$$
We call $A_i$ the $i$th {\it distance matrix} of $\Gamma.$
The matrix $A_1$ is often called  the {\it adjacency
matrix} of $\Gamma.$ We observe
(i) $A_0 = I$;
 (ii)
$\sum_{i=0}^D A_i = J$;
(iii)
$\overline{A_i} = A_i \;(0 \le i \le D)$;
(iv) $A_i^t = A_i  \;(0 \le i \le D)$;
(v) $A_iA_j = \sum_{h=0}^D p_{ij}^h A_h \;( 0 \le i,j \le D)
$,
where $I$ (resp. $J$) denotes the identity matrix 
(resp. all 1's matrix) in 
 $\hbox{Mat}_X(\C)$.
 Using these facts  we find
 $A_0,A_1,\ldots,A_D$
is a basis for a commutative subalgebra $M$ of 
$\mbox{Mat}_X(\C)$, called the 
{\it Bose-Mesner algebra} of $\Gamma$.
It turns out that $A_1$ generates $M$ \cite[p.~190]{bannai}.
By \cite[p.~45]{bcn}, $M$ has a second basis 
$E_0,E_1,\ldots,E_D$ such that
(i) $E_0 = |X|^{-1}J$;
(ii) $\sum_{i=0}^D E_i = I$;
(iii) $\overline{E_i} = E_i \;(0 \le i \le D)$;
(iv) $E_i^t =E_i  \;(0 \le i \le D)$;
(v) $E_iE_j =\delta_{ij}E_i  \;(0 \le i,j \le D)$.
We call $E_0, E_1, \ldots, E_D $  the {\it primitive idempotents}
of $\Gamma$.  

\medskip
\noindent
We  recall the eigenvalues
of  $\Gamma $.
Since $E_0,E_1,\ldots,E_D$ form a basis for  
$M$ there exist complex scalars $\theta_0,\theta_1,
\ldots,\theta_D$ such that
$A_1 = \sum_{i=0}^D \theta_iE_i$.
Observe
$A_1E_i = E_iA_1 =  \theta_iE_i$ for $0 \leq i \leq D$.
By \cite[p.~197]{bannai} the 
scalars $\theta_0,\theta_1,\ldots,\theta_D$ are
in $\R.$ Observe
$\theta_0,\theta_1,\ldots,\theta_D$ are mutually distinct 
since $A_1$ generates $M$. We call $\theta_i$  the {\it eigenvalue}
of $\Gamma$ associated with $E_i$ $(0 \leq i \leq D)$.
Observe 
\begin{eqnarray*}
V = E_0V+E_1V+ \cdots +E_DV \qquad \qquad {\rm (orthogonal\ direct\ sum}).
\end{eqnarray*}
For $0 \le i \le D$ the space $E_iV$ is the  eigenspace of $A_1$ associated 
with $\theta_i$.

\medskip
\noindent 
We now recall the Krein parameters.
Let $\circ $ denote the entrywise product in
$\mbox{Mat}_X(\C)$.
Observe
$A_i\circ A_j= \delta_{ij}A_i$ for $0 \leq i,j\leq D$,
so
$M$ is closed under
$\circ$. Thus there exist complex scalars
$q^h_{ij}$  $(0 \leq h,i,j\leq D)$ such
that
$$
E_i\circ E_j = |X|^{-1}\sum_{h=0}^D q^h_{ij}E_h
\qquad (0 \leq i,j\leq D).
$$
By \cite[p.~170]{Biggs}, 
$q^h_{ij}$ is real and nonnegative  for $0 \leq h,i,j\leq D$.
The $q^h_{ij}$ are called the {\it Krein parameters} of $\Gamma$.
The graph $\Gamma$ is said to be {\it $Q$-polynomial}
(with respect to the given ordering $E_0, E_1, \ldots, E_D$
of the primitive idempotents)
whenever for $0 \leq h,i,j\leq D$, 
$q^h_{ij}= 0$
(resp. 
$q^h_{ij}\not= 0$) whenever one of $h,i,j$ is greater than
(resp. equal to) the sum of the other two
\cite[p.~235]{bcn}. See
\cite{
wdw,
caugh1,
caugh2,
curtin3,
curtin4,
dickie1,
dickie2,
aap1}
 for background information on the $Q$-polynomial property.
For the rest of this section we assume $\Gamma$ is $Q$-polynomial
with respect to $E_0,E_1,\ldots,E_D$.

\medskip
\noindent
We  recall the dual Bose-Mesner algebra of $\Gamma.$
For the rest of this paper we fix
a vertex $x \in X.$ We view $x$ as a ``base vertex.''
For 
$ 0 \le i \le D$ let $E_i^*=E_i^*(x)$ denote the diagonal
matrix in $\hbox{Mat}_X(\C)$ with $(y,y)$-entry
\begin{equation}\label{DEFDEI}
{(E_i^*)_{yy} = \cases{1, & if $\partial(x,y)=i$\cr
0, & if $\partial(x,y) \ne i$\cr}} \qquad (y \in X).
\end{equation}
We call $E_i^*$ the  $i$th {\it dual idempotent} of $\Gamma$
 with respect to $x$ \cite[p.~378]{terwSub1}.
We observe
(i) $\sum_{i=0}^D E_i^*=I$;
(ii) $\overline{E_i^*} = E_i^*$ $(0 \le i \le D)$;
(iii) $E_i^{*t} = E_i^*$ $(0 \le i \le D)$;
(iv) $E_i^*E_j^* = \delta_{ij}E_i^* $ $(0 \le i,j \le D)$.
By these facts 
$E_0^*,E_1^*, \ldots, E_D^*$ form a 
basis for a commutative subalgebra
$M^*=M^*(x)$ of 
$\hbox{Mat}_X(\C).$ 
We call 
$M^*$ the {\it dual Bose-Mesner algebra} of
$\Gamma$ with respect to $x$ \cite[p.~378]{terwSub1}.
For $0 \leq i \leq D$ let $A^*_i = A^*_i(x)$ denote the diagonal
matrix in 
 $\hbox{Mat}_X(\C)$
with $(y,y)$-entry
$(A^*_i)_{yy}=\vert X \vert (E_i)_{xy}$ for $y \in X$.
Then $A^*_0, A^*_1, \ldots, A^*_D$ is a basis for $M^*$ 
\cite[p.~379]{terwSub1}.
Moreover
(i) $A^*_0 = I$;
(ii)
$\overline{A^*_i} = A^*_i \;(0 \le i \le D)$;
(iii) $A^{*t}_i = A^*_i  \;(0 \le i \le D)$;
(iv) $A^*_iA^*_j = \sum_{h=0}^D q_{ij}^h A^*_h \;( 0 \le i,j \le D)
$
\cite[p.~379]{terwSub1}.
We call 
 $A^*_0, A^*_1, \ldots, A^*_D$
the {\it dual distance matrices} of $\Gamma$ with respect to $x$.
The matrix
$A^*_1$ 
is often called  the {\it dual adjacency matrix} of $\Gamma$ with
respect to $x$.
The matrix $A^*_1$ generates $M^*$ \cite[Lemma 3.11]{terwSub1}.

\medskip
\noindent We recall the dual eigenvalues of $\Gamma$.
Since $E^*_0,E^*_1,\ldots,E^*_D$ form a basis for  
$M^*$ there exist complex scalars $\theta^*_0,\theta^*_1,
\ldots,\theta^*_D$ such that
$A^*_1 = \sum_{i=0}^D \theta^*_iE^*_i$.
Observe
$A^*_1E^*_i = E^*_iA^*_1 =  \theta^*_iE^*_i$ for $0 \leq i \leq D$.
By \cite[Lemma 3.11]{terwSub1} the 
scalars $\theta^*_0,\theta^*_1,\ldots,\theta^*_D$ are in $\R.$ 
The scalars $\theta^*_0,\theta^*_1,\ldots,\theta^*_D$ are mutually
distinct 
since $A^*_1$ generates $M^*$. We call $\theta^*_i$ the {\it dual eigenvalue}
of $\Gamma$ associated with $E^*_i$ $(0 \leq i\leq D)$.

\medskip
\noindent 
We recall the subconstituents of $\Gamma $.
From
(\ref{DEFDEI}) we find
\begin{equation}\label{DEIV}
E_i^*V = \mbox{span}\{\hat{y} \;|\; y \in X, \quad \partial(x,y)=i\}
\qquad (0 \le i \le D).
\end{equation}
By 
(\ref{DEIV})  and since
 $\{\hat{y}\;|\;y \in X\}$ is an orthonormal basis for $V$
 we find
\begin{eqnarray*}
\label{vsub}
V = E_0^*V+E_1^*V+ \cdots +E_D^*V \qquad \qquad 
{\rm (orthogonal\ direct\ sum}).
\end{eqnarray*}
For $0 \leq i \leq D$ the space $E^*_iV$ is the eigenspace
of $A^*_1$ associated with $\theta^*_i$.
We call $E_i^*V$ the $i$th {\it subconstituent} of $\Gamma$
with respect to $x$.

\medskip
\noindent
We recall the subconstituent algebra of $\Gamma $.
Let $T=T(x)$ denote the subalgebra of $\hbox{Mat}_X(\C)$ generated by 
$M$ and $M^*$. 
We call $T$ the {\it subconstituent algebra} 
(or {\it Terwilliger algebra}) of $\Gamma$ 
 with respect to $x$ \cite[Definition 3.3]{terwSub1}.
Observe that $T$ has finite dimension. Moreover $T$ is 
semisimple since it
is closed under the conjugate transponse map
\cite[p.~157]{CR}. By
\cite[Lemma 3.2]{terwSub1}
the following are relations in $T$:
\begin{eqnarray}
E^*_hA_iE^*_j&=&0 \quad \mbox{iff} \quad p^h_{ij}=0,
\qquad \qquad (0 \leq h,i,j \leq D),
\label{eq:triple2}
\\
E_hA^*_iE_j&=&0 \quad \mbox{iff} \quad q^h_{ij}=0, \qquad \qquad
(0 \leq h,i,j \leq D).
\label{eq:triple1}
\end{eqnarray}
See
\cite{curtin1,
curtin2,
curtin6,
egge1,
go,
go2,
hobart,
tanabe,
terwSub1,
terwSub2,
terwSub3}
for more information on the subconstituent
algebra.

\medskip
\noindent We recall the $T$-modules.
By a {\it T-module}
we mean a subspace $W \subseteq V$ such that $BW \subseteq W$
for all $B \in T.$ 
\noindent
Let $W$ denote a $T$-module and let 
$W'$ denote a  
$T$-module contained in $W$.
Then the orthogonal complement of $W'$ in $W$ is a $T$-module 
\cite[p.~802]{go2}.
It follows that each $T$-module
is an orthogonal direct sum of irreducible $T$-modules.
In particular $V$ is an orthogonal direct sum of irreducible $T$-modules.

\medskip
\noindent 
Let $W$ denote an irreducible $T$-module.
Observe that $W$ is the direct sum of the nonzero spaces among
$E^*_0W,\ldots, E^*_DW$. Similarly
$W$ is the direct sum 
 of the nonzero spaces among
$E_0W,\ldots,$ $ E_DW$.
By the {\it endpoint} of $W$ we mean
$\mbox{min}\lbrace i |0\leq i \leq D, \; E^*_iW\not=0\rbrace $.
By the {\it diameter} of $W$ we mean
$ |\lbrace i | 0 \leq i \leq D,\; E^*_iW\not=0 \rbrace |-1 $.
By the {\it dual endpoint} of $W$ we mean
$\mbox{min}\lbrace i |0\leq i \leq D, \; E_iW\not=0\rbrace $.
By
the {\it dual diameter} of $W$ we mean
$ |\lbrace i | 0 \leq i \leq D,\; E_iW\not=0 \rbrace |-1 $.
It turns out that the
diameter of $W$ is  equal to the dual diameter of
$W$
\cite[Corollary 3.3]{aap1}.
We finish this section with a comment.

\begin{lemma}
{\rm \cite[Lemma 3.4, Lemma 3.9, Lemma 3.12]{terwSub1}}
\label{lem:basic}
Let $W$ denote an irreducible $T$-module with endpoint $\rho$,
dual endpoint $\tau$, and diameter $d$.
Then $\rho,\tau,d$ are nonnegative integers such that $\rho+d\leq D$ and
$\tau+d\leq D$. Moreover the following (i)--(iv) hold.
\begin{enumerate}
\item 
$E^*_iW \not=0$ if and only if $\rho \leq i \leq \rho+d$, 
$ \quad (0 \leq i \leq D)$.
\item
$W = \sum_{h=0}^{d} E^*_{\rho+h}W \qquad (\mbox{orthogonal direct sum}). $
\item 
$E_iW \not=0$ if and only if $\tau \leq i \leq \tau+d$,
$ \quad (0 \leq i \leq D)$.
\item
$W = \sum_{h=0}^{d} E_{\tau+h}W \qquad (\mbox{orthogonal direct sum}). $
\end{enumerate}
\end{lemma}

\section{A restriction on the intersection numbers}

From now on we impose the following restriction on the
intersection numbers of $\Gamma$.

\begin{assumption}
\label{def:sdcp}
\rm
We fix 
$b,\beta \in \C$ such that $b\not=1$, and
assume $\Gamma$ has classical parameters $(D,b,\alpha,\beta)$
with $\alpha=b-1$.
 This means that the intersection numbers
of $\Gamma$ satisfy
\begin{eqnarray*}
c_i &=& b^{i-1}\frac{b^i-1}{b-1}, \\
b_i &=& (\beta+1-b^i)\frac{b^D-b^i}{b-1}
\end{eqnarray*}
for $0 \leq i \leq D$
\cite[p.~193]{bcn}.
We remark that 
$b$ is an integer and $b\not=0$, $b\not=-1$ \cite[Proposition~6.2.1]{bcn}.
For notational convenience we fix $q \in \C$ such that
\begin{eqnarray*}
b = q^2.
\end{eqnarray*}
We note that $q$ is nonzero and not a root of unity.
\end{assumption}

\begin{remark}
\label{rem:qp}
\rm Referring to Assumption \ref{def:sdcp}, the restriction 
$\alpha=b-1$
 implies that
$\Gamma$ is formally
self-dual \cite[Corollary~8.4.4]{bcn}.
Consequently there exists an ordering
$E_0, E_1, \ldots, E_D$ of the primitive idempotents of
$\Gamma$,
with respect to which
the Krein parameter $q^h_{ij}$ is
equal to the intersection number $p^h_{ij}$ for
$0 \leq h,i,j\leq D$.
In particular $\Gamma$ is $Q$-polynomial with respect to
$E_0, E_1, \ldots, E_D$.
We fix this ordering of the primitive idempotents for the rest of
the paper. 
\end{remark}

\begin{remark}
\rm
In the notation of Bannai and Ito \cite[p. 263]{bannai},
the $Q$-polynomial structure from
Remark \ref{rem:qp}
is type I with $s=0, s^*=0$.
\end{remark}

\begin{example}
\label{ex:nicegraph}
\rm 
The following distance-regular graphs satisfy
 Assumption
\ref{def:sdcp}:
the
bilinear forms graph \cite[p.~280]{bcn},
the alternating forms graph \cite[p.~282]{bcn},
the Hermitean forms graph \cite[p.~285]{bcn},
the
quadratic forms graph \cite[p.~290]{bcn},
the affine $E_6$ graph \cite[p.~340]{bcn},
and the extended ternary Golay code graph \cite[p.~359]{bcn}.
\end{example}

\noindent
With reference to Assumption \ref{def:sdcp}
we will display an action of $\boxtimes_q$ on the
standard module of $\Gamma$. To describe this
action we define eight matrices in 
 $\hbox{Mat}_X(\C)$, called
\begin{eqnarray}
\label{eq:list}
A,\quad
A^*,\quad
B,\quad
B^*, \quad
K,\quad
K^*,\quad
\Phi,\quad
\Psi.
\end{eqnarray}
These matrices will be defined in the next two sections.

\section{The matrices $A$ and $A^*$}

\noindent In this section we define the matrices $A, A^*$ and
discuss their properties.
We start with a comment.

\begin{lemma}
{\rm \cite[Corollary~8.4.4]{bcn}}
\label{lem:alpha}
With reference to Assumption \ref{def:sdcp}, there
exist $\alpha_0,\alpha_1 \in \C$ such that 
each of $\theta_i$, $\theta^*_i$ is 
 $\alpha_0+\alpha_1 q^{D-2i}$ for $0 \leq i \leq D$.
Moreover $\alpha_1\not=0$.
\end{lemma}

\begin{definition}
\label{def:aas} \rm
With reference to Assumption \ref{def:sdcp}
we define $A,A^* \in 
\hbox{Mat}_X(\C)$ so that
\begin{eqnarray*}
A_1&=&\alpha_0I+\alpha_1A,
\\
A^*_1&=&\alpha_0I+\alpha_1A^*,
\end{eqnarray*} 
where $\alpha_0,\alpha_1$ are from
Lemma \ref{lem:alpha}.
Thus for $0 \leq i \leq D$ the 
space $E_iV$ (resp. $E^*_iV$) is an
eigenspace of $A$ (resp. $A^*$) with eigenvalue
 $q^{D-2i}$.
\end{definition}

\begin{lemma}
\label{lem:td}
With reference to Assumption \ref{def:sdcp}
and Definition 
\ref{def:aas},
 the following (i), (ii) hold for all $0 \leq i,j\leq D$
such that $|i-j|>1$:
\begin{enumerate}
\item $E^*_iAE^*_j = 0$,
\item $E_iA^*E_j = 0$. 
\end{enumerate}
\end{lemma}
\noindent {\it Proof:} (i) 
We have $p^i_{1j}=0$ since $ \vert i-j\vert >1$, so
$E^*_iA_1E^*_j=0$ in view of
(\ref{eq:triple2}). The result now follows using
the first equation of Definition
\ref{def:aas}.
\\
\noindent (ii) Similar to the proof of (i) above.
\hfill $\Box $ \\

\noindent The following is essentially a special case of
\cite[Lemma 5.4]{terwSub3}.

\begin{lemma} {\rm
\cite[Lemma 5.4]{terwSub3}}
With reference to Assumption \ref{def:sdcp} and
Definition \ref{def:aas} the matrices $A,A^*$ satisfy the
$q$-Serre relations
\begin{eqnarray}\label{eq:qs1}
A^3A^*-
\lbrack 3\rbrack_q A^2A^*A
+
\lbrack 3\rbrack_q AA^*A^2
-A^*A^3 &=&0,\\
\label{eq:qs2}
A^{*3}A-
\lbrack 3\rbrack_q A^{*2}AA^*
+
\lbrack 3\rbrack_q A^*AA^{*2}
-AA^{*3} &=&0.
\end{eqnarray}
\end{lemma}
\noindent {\it Proof:}
We first show (\ref{eq:qs1}). By the last sentence in Definition
\ref{def:aas}, for $0 \leq i \leq D$ we have
$AE_i=E_iA=\sigma_iE_i$ where $\sigma_i=q^{D-2i}$.
Let $C$ denote the expression on the left in (\ref{eq:qs1}).
We show $C=0$. Since $I=E_0+\cdots +E_D$ it suffices to
show $E_iCE_j=0$ for $0 \leq i,j\leq D$. Let $i,j$ be given.
By our preliminary comment and the definition of $C$ we find
$E_iCE_j=E_iA^*E_j\alpha_{ij}$ where
\begin{eqnarray}
\alpha_{ij} &=& 
\sigma_i^3
-
\lbrack 3\rbrack_q \sigma^2_i \sigma_j
+
\lbrack 3\rbrack_q \sigma_i\sigma^2_j
-\sigma^3_j 
\nonumber
\\
&=& (\sigma_i-\sigma_jq^2)(\sigma_i-\sigma_j)(\sigma_i-\sigma_jq^{-2}).
\label{eq:tprod}
\end{eqnarray}
If $|i-j|>1$ then $E_iA^*E_j=0$ by Lemma
\ref{lem:td}(ii).
If $|i-j|\leq 1$ then $\alpha_{ij}=0$ by
(\ref{eq:tprod}) and the definition of $\sigma_0,\ldots,\sigma_D$.
In either case $E_iCE_j=0$ as desired.
It follows
that $C=0$ and line
(\ref{eq:qs1}) is proved.
The proof of (\ref{eq:qs2}) is similar to the proof of
(\ref{eq:qs1}).
\hfill $\Box $ \\

\noindent We finish this section with a comment.
\begin{lemma}
\label{lem:aasgen}
With reference to Assumption \ref{def:sdcp} and Definition
\ref{def:aas} the matrices $A,A^*$ together generate
$T$.
\end{lemma}
\noindent {\it Proof:}
By definition $T$ is generated by $M$ and $M^*$.
The algebra $M$ (resp. $M^*$) is generated by
$A_1$ (resp. $A^*_1$)
and hence by $A$ (resp. $A^*$) 
in view of Definition
\ref{def:aas}.
The result follows.
\hfill $\Box $ \\

\section{The matrices $B$, $B^*$, $K$, $K^*$, $\Phi$, $\Psi$}

\noindent In the previous section we defined the matrices
$A, A^*$. In this section we define the remaining matrices 
from the list
(\ref{eq:list}).

\begin{definition}
\label{def:updown}
\rm
With reference to Assumption \ref{def:sdcp}, for $-1\leq i,j\leq D$
we define
\begin{eqnarray*}
V_{i,j}^{\downarrow \downarrow} &=& 
(E^*_0V+\cdots+E^*_iV)\cap (E_0V+\cdots+E_jV),
\\
V_{i,j}^{\uparrow \downarrow} &=& 
(E^*_DV+\cdots+E^*_{D-i}V)\cap (E_0V+\cdots+E_jV),
\\
V_{i,j}^{\downarrow \uparrow} &=& 
(E^*_0V+\cdots+E^*_iV)\cap (E_DV+\cdots+E_{D-j}V),
\\
V_{i,j}^{\uparrow \uparrow} &=& 
(E^*_DV+\cdots+E^*_{D-i}V)\cap (E_DV+\cdots+E_{D-j}V).
\end{eqnarray*}
In each of the above four equations we interpret the
right-hand side to be 0 if $i=-1$ or $j=-1$. 
\end{definition}

\begin{definition}
\label{def:vtilde}
\rm
With reference to Assumption \ref{def:sdcp} and Definition
\ref{def:updown}, for $\eta, \mu \in \lbrace \downarrow,
\uparrow \rbrace$ and  $0 \leq i,j\leq D$ we have
$
V^{\eta \mu}_{i-1,j} \subseteq V^{\eta \mu}_{i,j}$
and
$
V^{\eta \mu}_{i,j-1}  \subseteq V^{\eta \mu}_{i,j}
$. Therefore
\begin{eqnarray*}
V^{\eta \mu}_{i-1,j}+
V^{\eta \mu}_{i,j-1} \subseteq V^{\eta \mu}_{i,j}.
\end{eqnarray*}
Referring to the above inclusion, we define ${\tilde V}^{\eta \mu}_{i,j}$
 to be the orthogonal complement of the left-hand side in the
 right-hand side; that
is
\begin{eqnarray*}
 {\tilde V}^{\eta \mu}_{i,j}=(
V^{\eta \mu}_{i-1,j}+
V^{\eta \mu}_{i,j-1})^\perp \cap V^{\eta \mu}_{i,j}.
\end{eqnarray*}
\end{definition}

\noindent The following result is a mild generalization of
\cite[Corollary 5.8]{ds}.

\begin{lemma}
With reference to Assumption \ref{def:sdcp} and Definition
\ref{def:vtilde} the following holds 
 for $\eta,\mu \in \lbrace \downarrow, \uparrow \rbrace$:
\begin{eqnarray*}
V = \sum_{i=0}^D \sum_{j=0}^D {\tilde V}^{\eta \mu}_{i,j}
\qquad \qquad (\mbox{direct sum}).
\end{eqnarray*}
\end{lemma}
\noindent {\it Proof:}
For $\eta=\downarrow $, $ \mu=\downarrow $
 this is just \cite[Corollary~5.8]{ds}.
For general values of $\eta,\mu$, in the proof
of \cite[Corollary~5.8]{ds} replace the sequence
$E^*_0,\ldots,E^*_D$ 
(resp. $E_0,\ldots,E_D$)
by $E^*_D,\ldots,E^*_0$ (resp. 
$E_D,\ldots,E_0$)
 if $\eta=\uparrow$
(resp. $\mu=\uparrow$).
\hfill $\Box $ \\

\begin{definition}
\label{def:listdef}
\rm
With reference to Assumption \ref{def:sdcp} and
Definition \ref{def:vtilde},
we define $B$, $B^*$, $K$, $K^*$, $\Phi$, $\Psi$ to be 
the unique matrices
in
 $\hbox{Mat}_X(\C)$
that satisfy the requirements of the following table
for $0 \leq i,j\leq D$.

\medskip
\centerline{
\begin{tabular}[t]{c|c}
        {\rm The matrix}
     &
  {\rm is $0$ on}
\\  	\hline 
       	 $B-q^{i-j}I$
         &
  ${\tilde V}_{i,j}^{\downarrow \uparrow}$
\\
 $B^*-q^{j-i}I$  
         &
 ${\tilde V}_{i,j}^{\uparrow \downarrow}$
\\
 $K-q^{i-j}I$  
         &
 ${\tilde V}_{i,j}^{\downarrow \downarrow}$
\\
 $K^*-q^{i-j}I$  
         & 
 ${\tilde V}_{i,j}^{\uparrow \uparrow}$
\\
$\Phi-q^{i+j-D}I$  
         &
 ${\tilde V}_{i,j}^{\downarrow \downarrow}$
\\
 $\Psi-q^{i+j-D}I$ 
&
 ${\tilde V}_{i,j}^{\downarrow \uparrow}$
\end{tabular}}

\medskip
\end{definition}

\section{An action of $\boxtimes_q$ on the standard module of $\Gamma$}

We now state our main result, in which we
display an action of $\boxtimes_q$ on the standard module $V$
of $\Gamma$.

\begin{theorem}
\label{thm:mres}
With reference to Assumption \ref{def:sdcp}, there exists 
a $\boxtimes_q$-module structure on $V$ such
that the generators $x_{ij}$
act as follows:

\medskip
\centerline{
\begin{tabular}[t]{c|cccccccc}
        {\rm generator}  
       	& $x_{01}$  
         & $x_{12}$  
         & $x_{23}$  
         & $x_{30}$  
         & $x_{02}$  
         & $x_{13}$    
	\\
	\hline 
{\rm action on $V$} 
&  $A\Phi \Psi^{-1}$ & $B \Phi^{-1}$ & $A^* \Phi\Psi$ & $B^*\Phi^{-1}$ & $K\Psi^{-1}$ & $K^*\Psi$ 
\end{tabular}}
\end{theorem}

\noindent The proof of Theorem
\ref{thm:mres}
is given at the end of this section.
First we need some lemmas.

\begin{lemma}
\label{lem:localtet}
With reference to Assumption \ref{def:sdcp}, let $W$
denote an irreducible $T$-module with endpoint $\rho$,
dual endpoint $\tau$, and diameter $d$. Then there exists
a unique $\boxtimes_q$-module structure on $W$ such that
the generators $x_{01}$, $x_{23}$ act as
$Aq^{d-D+2\tau}$, $A^*q^{d-D+2\rho}$ respectively.
This $\boxtimes_q$-module structure is irreducible and
type $1$.
\end{lemma}
\noindent {\it Proof:}
The matrices $A,A^*$ satisfy the $q$-Serre relations
(\ref{eq:qs1}),
(\ref{eq:qs2}).   
These relations are homogeneous so they
still hold if $A,A^*$ are replaced by
$Aq^{d-D+2\tau}$, $A^*q^{d-D+2\rho}$ respectively.
Therefore there exists an ${\mathcal A}_q$-module
structure on $W$ such that the standard generators act as
$Aq^{d-D+2\tau}$ and  $A^*q^{d-D+2\rho}$.
The ${\mathcal A}_q$-module $W$ is irreducible since
$A,A^*$ generate $T$ and since the $T$-module $W$
is irreducible. By Lemma
\ref{lem:basic}(iii),(iv)
the action of $A$ on $W$ is semisimple with eigenvalues
$q^{D-2\tau-2i}$ $(0 \leq i \leq d)$. Therefore the
action of 
$Aq^{d-D+2\tau}$ on $W$ is semisimple with eigenvalues
$q^{d-2i}$ $(0 \leq i \leq d)$.
 By Lemma
\ref{lem:basic}(i),(ii)
the action of $A^*$ on $W$ is semisimple with eigenvalues
$q^{D-2\rho-2i}$ $(0 \leq i \leq d)$. Therefore the
action of 
$A^*q^{d-D+2\rho}$ on $W$ is semisimple with eigenvalues
$q^{d-2i}$ $(0 \leq i \leq d)$. By these comments
and the first paragraph of Section 6
the ${\mathcal A}_q$-module $W$ is NonNil and type $(1,1)$.
So far we have shown that the ${\mathcal A}_q$-module
$W$ is irreducible, NonNil, and type $(1,1)$.
Combining this with Theorem
\ref{thm:1} we obtain the result.
\hfill $\Box $ \\

\begin{lemma}
\label{lem:xijeig}
With reference to Assumption \ref{def:sdcp}, let $W$
denote an irreducible $T$-module with endpoint $\rho$,
dual endpoint $\tau$, and diameter $d$. Consider the
$\boxtimes_q$-module structure on $W$ from Lemma
\ref{lem:localtet}. For each generator $x_{rs}$ of $\boxtimes_q$
and for $0 \leq i \leq d$, the eigenspace of $x_{rs}$ on $W$
associated with the eigenvalue $q^{d-2i}$ is given in 
the following table.

\medskip
\centerline{
\begin{tabular}[t]{c|c|c}
        $r$ & $s$ & {\rm eigenspace of $x_{rs}$ for the eigenvalue $q^{d-2i}$}
\\    
\hline
$0$ & $1$ & $E_{\tau+i}W$
\\
$1$ & $2$ & $(E^*_{\rho}W+\cdots+E^*_{\rho+d-i}W)\cap
(E_{\tau+d-i}W+\cdots+E_{\tau+d}W)$
\\
$2$ & $3$ & $E^*_{\rho+i}W$
\\
$3$ & $0$ & $(E^*_{\rho+d-i}W+\cdots+E^*_{\rho+d}W)\cap
(E_{\tau}W+\cdots+E_{\tau+d-i}W)$
\\
$0$ & $2$ & $(E^*_{\rho}W+\cdots+E^*_{\rho+d-i}W)\cap
(E_{\tau}W+\cdots+E_{\tau+i}W)$
\\
$1$ & $3$ & $(E^*_{\rho+i}W+\cdots+E^*_{\rho+d}W)\cap
(E_{\tau+d-i}W+\cdots+E_{\tau+d}W)$
\end{tabular}}

\end{lemma}
\noindent {\it Proof:}
Referring to the table, we first verify row $(r,s)=(0,1)$.
By Lemma
\ref{lem:localtet}
the generator $x_{01}$ acts on $W$ as $Aq^{d-D+2\tau}$.
By Lemma 
\ref{lem:basic}(iii),(iv)
the space $E_{\tau+i}W$ is the eigenspace
of $A$ on $W$ for the eigenvalue $q^{D-2\tau-2i}$. By
these comments $E_{\tau+i}W$ is the eigenspace of $x_{01}$ on
$W$ for the eigenvalue $q^{d-2i}$. We have now verified
row $(r,s)=(0,1)$.
Next we verify row
$(r,s)=(2,3)$.
By Lemma
\ref{lem:localtet}
the generator $x_{23}$ acts on $W$ as $A^*q^{d-D+2\rho}$.
By Lemma 
\ref{lem:basic}(i),(ii)
the space $E^*_{\rho+i}W$ is the eigenspace
of $A^*$ on $W$ for the eigenvalue $q^{D-2\rho-2i}$. By
these comments $E^*_{\rho+i}W$ is the eigenspace of $x_{23}$ on
$W$ for the eigenvalue $q^{d-2i}$. We have now verified
row $(r,s)=(2,3)$.
The remaining rows are valid by
\cite[Theorem 16.4]{qtet}.
\hfill $\Box $ \\

\noindent The following result is a mild generalization of
\cite[Lemma 6.1]{ds}.

\begin{lemma}
\label{lem:xijw}
With reference to Assumption \ref{def:sdcp}, let $W$
denote an irreducible $T$-module with endpoint $\rho$,
dual endpoint $\tau$, and diameter $d$.
Then the following (i)--(iv) hold for $0 \leq i \leq d$.
\begin{enumerate}
\item The space
\begin{eqnarray*}
(E^*_{\rho}W+\cdots+E^*_{\rho+d-i}W)\cap
(E_{\tau+d-i}W+\cdots+E_{\tau+d}W)
\end{eqnarray*}
is contained in ${\tilde V}^{\downarrow \uparrow}_{\rho+d-i,D-d-\tau+i}$.
\item The space
\begin{eqnarray*}
(E^*_{\rho+d-i}W+\cdots+E^*_{\rho+d}W)\cap
(E_{\tau}W+\cdots+E_{\tau+d-i}W)
\end{eqnarray*}
is contained in 
 ${\tilde V}^{\uparrow \downarrow}_{D-d-\rho+i,\tau+d-i}$.
\item The space
\begin{eqnarray*}
(E^*_{\rho}W+\cdots+E^*_{\rho+d-i}W)\cap
(E_{\tau}W+\cdots+E_{\tau+i}W)
\end{eqnarray*}
is contained in
 ${\tilde V}^{\downarrow \downarrow}_{\rho+d-i,\tau+i}$.
\item The space
\begin{eqnarray*}
(E^*_{\rho+i}W+\cdots+E^*_{\rho+d}W)\cap
(E_{\tau+d-i}W+\cdots+E_{\tau+d}W)
\end{eqnarray*}
is contained in 
 ${\tilde V}^{\uparrow \uparrow}_{D-\rho-i,D-d-\tau+i}$.
\end{enumerate}
\end{lemma}
\noindent {\it Proof:}
Assertion 
(iii) is just \cite[Lemma 6.1]{ds}.
To get  (i), in the proof of
\cite[Lemma 6.1]{ds} replace
the sequence $E_0,\ldots, E_D$ by $E_D,\ldots,E_0$.
To  get  (ii), in the proof of
\cite[Lemma 6.1]{ds} replace
 $E^*_0,\ldots, E^*_D$ by $E^*_D,\ldots,E^*_0$.
To get (iv), in the proof of
\cite[Lemma 6.1]{ds} replace
 $E^*_0,\ldots, E^*_D$ 
(resp. $E_0,\ldots, E_D$) 
by $E^*_D,\ldots,E^*_0$
(resp. $E_D,\ldots,E_0$).
\hfill $\Box $ \\

\begin{lemma}
\label{lem:abx}
With reference to Assumption \ref{def:sdcp}, let $W$
denote an irreducible $T$-module with endpoint $\rho$,
dual endpoint $\tau$, and diameter $d$. Consider the
$\boxtimes_q$-module structure on $W$ from Lemma
\ref{lem:localtet}. 
In the table below, each row contains a matrix in
$\hbox{Mat}_X(\C)$ 
 and
an element of $\boxtimes_q$. The action of these two objects
on $W$ coincide.

\medskip
\centerline{
\begin{tabular}[t]{c|c}
        {\rm matrix}
     &
  {\rm element of $\boxtimes_q$} 
\\	\hline 
       	 $A$  
         &
  $q^{D-d-2\tau}x_{01}$
\\
 $B$  
         &
 $q^{d-D+\rho+\tau}x_{12}$
\\
 $A^*$  
         &
 $q^{D-d-2\rho}x_{23}$
\\
 $B^*$  
         & 
 $q^{d-D+\rho+\tau}x_{30}$
\\
$K$  
         &
 $q^{\rho-\tau}x_{02}$
\\
 $K^*$ 
&
 $q^{\tau-\rho}x_{13}$
\\
$\Phi$  
         &
 $q^{d-D+\rho+\tau}1$
\\
 $\Psi$ 
&
 $q^{\rho-\tau}1$
\end{tabular}}

\medskip
\end{lemma}
\noindent {\it Proof:}
By Lemma 
\ref{lem:localtet} the expressions
$A-q^{D-d-2\tau}x_{01}$
and 
$A^*-q^{D-d-2\rho}x_{23}$
are each 0 on $W$.
Next we show that
$B-q^{d-D+\rho+\tau}x_{12}$ is 0 on $W$.
To this end we pick $w \in W$ and show
$Bw=q^{d-D+\rho+\tau}x_{12}w$.
Recall that $x_{12}$ is semisimple on $W$ with eigenvalues
$q^{d-2i}$ $(0 \leq i \leq d)$.
Therefore without loss of generality we may assume that 
there exists an integer $i$ $(0 \leq i\leq d)$ such 
that $x_{12}w=q^{d-2i}w$. By row $(r,s)=(1,2)$ in the table
of Lemma
\ref{lem:xijeig} and by Lemma
\ref{lem:xijw}(i), we find
$w \in {\tilde V}^{\downarrow \uparrow}_{\rho+d-i,D-d-\tau+i}$.
By this and the first row in the table of
Definition 
\ref{def:listdef}
we find
$Bw=q^{2d-D+\rho+\tau-2i}w$. From these comments
we find
$Bw=q^{d-D+\rho+\tau}x_{12}w$ as desired. We have now shown
that $B-q^{d-D+\rho+\tau}x_{12}$ is 0 on $W$.
Similarly one shows that each of 
 $B^*-q^{d-D+\rho+\tau}x_{30}$,
 $K-q^{\rho-\tau}x_{02}$,
 $K^*-q^{\tau-\rho}x_{13}$ is 0 on $W$.
We now show that
$\Phi-q^{d-D+\rho+\tau}I$ is 0 on $W$.
To this end we pick $v \in W$ and show
$\Phi v=q^{d-D+\rho+\tau}v$.
Recall that $x_{02}$ is semisimple on $W$ with eigenvalues
$q^{d-2i}$ $(0 \leq i \leq d)$.
Therefore without loss of generality we may assume
that there exists an integer $i$ $(0 \leq i\leq d)$ such 
that $x_{02}v=q^{d-2i}v$. By row $(r,s)=(0,2)$ in the table
of Lemma
\ref{lem:xijeig} and by Lemma
\ref{lem:xijw}(iii), we find
$v \in {\tilde V}^{\downarrow \downarrow}_{\rho+d-i,\tau+i}$.
By this and the second to the last row in the table of
Definition \ref{def:listdef} we find
$\Phi v=q^{d-D+\rho+\tau}v$ as desired.
 We have now shown
that $\Phi-q^{d-D+\rho+\tau}I$ is 0 on $W$.
Similarly one shows that  
 $\Psi-q^{\rho-\tau}I$ is 0 on $W$.
\hfill $\Box $ \\

\begin{corollary}
\label{cor:abx}
With reference to Assumption \ref{def:sdcp}, let $W$
denote an irreducible $T$-module and 
consider the
$\boxtimes_q$-action on $W$ from Lemma
\ref{lem:localtet}. 
In the table below, each column  contains a
generator for 
$\boxtimes_q$ and a matrix in
$\hbox{Mat}_X(\C)$.
The action of these two objects
on $W$ coincide.

\medskip
\centerline{
\begin{tabular}[t]{c|cccccccc}
        {\rm generator}  
       	& $x_{01}$  
         & $x_{12}$  
         & $x_{23}$  
         & $x_{30}$  
         & $x_{02}$  
         & $x_{13}$    
	\\
	\hline 
{\rm matrix} 
&  $A\Phi \Psi^{-1}$ & $B \Phi^{-1}$ & $A^* \Phi\Psi$ & $B^*\Phi^{-1}$ & $K\Psi^{-1}$ & $K^*\Psi$ 
\end{tabular}}
\end{corollary}

\noindent {\it Proof:}
Immediate from Lemma
\ref{lem:abx}.
\hfill $\Box $ \\

\noindent It is now a simple matter to prove Theorem
\ref{thm:mres}.

\medskip
\noindent {\it Proof of Theorem
\ref{thm:mres}:}
The standard module $V$ decomposes into a direct sum of
irreducible $T$-modules. Each irreducible $T$-module
in this decomposition supports a $\boxtimes_q$-module
structure from Lemma
\ref{lem:localtet}.
Combining these $\boxtimes_q$-modules we get a $\boxtimes_q$-module
structure on $V$. 
It remains to show that this $\boxtimes_q$-module satisfies the
requirements of 
Theorem \ref{thm:mres}.
This is the case since by Corollary 
\ref{cor:abx},
for each column in
 the table of Theorem \ref{thm:mres} the 
given $\boxtimes_q$ generator and the matrix beneath it
coincide
on each of the irreducible $T$-modules in the
above decomposition and hence on $V$.
\hfill $\Box $ \\

\begin{remark}
\label{rem:cp}
\rm In Theorem 
\ref{thm:mres} we displayed an action of 
$\boxtimes_q$ on the standard module $V$
of $\Gamma$.
In Proposition
\ref{prop:uqhom} we displayed four $\C$-algebra homomorphisms
from 
$U_q({\widehat{\mathfrak{sl}}_2})$ 
to $\boxtimes_q$. Using these homomorphisms to pull back the 
$\boxtimes_q$-action we obtain four 
$U_q({\widehat{\mathfrak{sl}}_2})$-module
structures on $V$.
\end{remark}

\section{How $\boxtimes_q$ is related to $T$}

\noindent
In Theorem 
\ref{thm:mres} we displayed an action of 
$\boxtimes_q$ on the standard module
of $\Gamma$; observe that this action induces
a $\C$-algebra homomorphism 
$\boxtimes_q \to 
\hbox{Mat}_X(\C)$ which we will denote by $\vartheta$.
In this section
we clarify how the
image $\vartheta(\boxtimes_q)$
is related to 
the subconstituent
algebra $T$.

\begin{lemma}
\label{lem:central}
With reference to Assumption \ref{def:sdcp}, the following
(i), (ii) hold.
\begin{enumerate}
\item
Each of the matrices from the list
(\ref{eq:list}) is contained in $T$.
\item
Each of $\Phi, \Psi$ is contained in the center $Z(T)$.
\end{enumerate}
\end{lemma}
\noindent {\it Proof:} 
(i) By Lemma \ref{lem:abx}
each matrix in the list
(\ref{eq:list}) leaves invariant every irreducible $T$-module.
Let 
$T'$ denote the set of matrices in
$\hbox{Mat}_X(\C)$ that leave invariant every irreducible
$T$-module. We observe that $T'$ is a subalgebra of
$\hbox{Mat}_X(\C)$ that contains $T$ as well as each
matrix in the list
(\ref{eq:list}).
We show that $T=T'$.
To this end we first show that $T'$ is semisimple.
By the construction each irreducible $T$-module is
an irreducible $T'$-module. We mentioned in Section 7 that
the standard module
$V$ is a direct sum of irreducible $T$-modules.
Therefore $V$ is a direct sum of irreducible
$T'$-modules, so $T'$ is semisimple. 
Next, let $W_1$, $W_2$ denote irreducible $T$-modules.
We claim that any isomorphism of $T$-modules
$\gamma: W_1 \to W_2$ is an isomorphism of $T'$-modules.
This is readily checked using the fact that
$\lbrace w+\gamma(w) \;|\;w \in W_1 \rbrace$ 
is an irreducible $T$-module and therefore invariant
under $T'$.
By our above comments the vector spaces 
$T$ and $T'$ have
the same dimension; this dimension is
$\sum_{\lambda} {d_{\lambda}^2}$ where the sum is over all
isomorphism classes $\lambda $ of irreducible $T$-modules
and $d_{\lambda}$ denotes the dimension of an irreducible
$T$-module in the isomorphism class $\lambda$.
Since $T'$ contains $T$ and they have the same
dimension we find $T=T'$. The result follows.
\\
\noindent (ii) By Lemma
\ref{lem:abx} 
each of $\Phi, \Psi$ acts as a scalar multiple
of the identity on every irreducible $T$-module.
\hfill $\Box $

\begin{theorem}
With reference to Assumption \ref{def:sdcp}
the following (i), (ii) hold.
\begin{enumerate}
\item
The image $\vartheta(\boxtimes_q)$ is contained in $T$.
\item
$T$ is generated by 
 $\vartheta(\boxtimes_q)$ together with $\Phi, \Psi$.
\end{enumerate}
\end{theorem}
\noindent {\it Proof:} 
Combine 
Lemma
\ref{lem:aasgen},
Theorem \ref{thm:mres},
 and
Lemma \ref{lem:central}.
\hfill $\Box $

\section{Directions for further research}

\noindent In this section we give some suggestions 
for further research.

\begin{problem}
\rm
For 
the spaces in Definition
\ref{def:updown}, find a combinatorial interpretation and 
an attractive basis.
\end{problem}

\begin{problem}
\rm
With reference to Assumption \ref{def:sdcp},
the matrices $\Phi, \Psi$ commute by 
Lemma \ref{lem:central}(ii) and they
are 
semisimple by Definition
\ref{def:listdef}.
 Therefore the
standard module of 
$\Gamma$ decomposes into a direct sum
of their common eigenspaces.
For these common eigenspaces find
a combinatorial interpretation 
and an attractive basis.
\end{problem}

\begin{problem}
\rm
With reference to Assumption \ref{def:sdcp},
for $y,z \in X$ and for each of
$B$, $B^*$, $K$, $K^*$, $\Phi$, $\Psi$ find the $(y,z)$-entry
in terms of
the distances $\partial(x,y)$,
$\partial(y,z)$, $\partial(z,x)$ ($x = \mbox{base vertex from Section 7}$)
 and
other combinatorial parameters as needed.
When is this entry 0?
\end{problem}

\begin{problem} \rm
Find all the distance-regular graphs 
that have classical parameters $(D,b,\alpha,\beta)$ and
$b\not=1$, $\alpha=b-1$.
 Some examples
are given in 
 Example
\ref{ex:nicegraph}.
\end{problem}

\begin{problem} \rm 
The finite-dimensional irreducible 
$U_q({\widehat{\mathfrak{sl}}_2})$-modules 
are classified by V. Chari and A. Pressley
\cite{charp}; see also
\cite{ding}, 
\cite{tarasov}.
Use this and Remark 
\ref{rem:cp}
to describe the irreducible
$T$-modules for each of the 
graphs in Example
\ref{ex:nicegraph}.
\end{problem}

\begin{conjecture}
With reference to Assumption \ref{def:sdcp}, 
for $0 \leq i,j\leq D$ the spaces 
${\tilde V}^{\downarrow \downarrow}_{ij}$
and 
${\tilde V}^{\uparrow \uparrow}_{rs}$ are
orthogonal unless $i+r=D$
and $j+s=D$.
Moreover 
${\tilde V}^{\downarrow \uparrow}_{ij}$
and 
${\tilde V}^{\uparrow \downarrow}_{rs}$ are
orthogonal unless $i+r=D$
and $j+s=D$.
\end{conjecture}

\begin{problem}
\rm
With reference to Assumption \ref{def:sdcp}, 
note by Lemma
\ref{lem:abx} that 
the following are equivalent:
(i) 
for each irreducible $T$-module the
endpoint and dual endpoint coincide;
(ii) $\Psi=I$.
For which of the graphs in
Example \ref{ex:nicegraph} do these equivalent
conditions hold?
\end{problem}

\begin{conjecture}
\label{lem:transpose}
With reference to Assumption \ref{def:sdcp}, 
 each of $\Phi,\Psi$ is symmetric and
\begin{eqnarray*}
B^t=B^*, \qquad \qquad K^t=K^{*-1}.
\end{eqnarray*}
\end{conjecture}

\noindent Under Assumption  
\ref{def:sdcp} we displayed an action of $\boxtimes_q$
on the standard module of $\Gamma$.
For the moment replace Assumption 
\ref{def:sdcp}
by the weaker assumption that 
$\Gamma$ is $Q$-polynomial. We suspect 
that there is still a natural action of $\boxtimes_q$
(or 
$U_q({\widehat{\mathfrak{sl}}_2})$, 
$U_q(\mathfrak{sl}_2)$,
${\widehat{\mathfrak{sl}}_2}$, 
$\mathfrak{sl}_2,
\ldots $ 
in degenerate cases) on the standard module of $\Gamma$.
It is premature for us to guess how
this action behaves
in every case, but the general idea is conveyed in
the following two conjectures.

\begin{conjecture}
\label{conj:cpar}
Assume $\Gamma$ has classical parameters
$(D,b,\alpha,\beta)$ with $b\not=1$.
In order to avoid degenerate situations,
assume that $\Gamma$ is not a dual polar graph 
\cite[p.~274]{bcn}. 
Then for $b=q^2$ there exists a $\boxtimes_q$-action
on the standard module of $\Gamma$
for which 
the adjacency matrix acts as 
a $Z(T)$-linear combination of 
$1, x_{01}, x_{12}$
and the dual adjacency matrix acts as  
a $Z(T)$-linear combination of 
$1, x_{23}$. We recall that $Z(T)$ denotes the center of $T$.
\end{conjecture}

\begin{conjecture}
Assume $\Gamma$ is $Q$-polynomial, with
eigenvalues
$\theta_i$ and dual eigenvalues $\theta^*_i$.
Recall that
the expressions
\begin{eqnarray*}
\frac{\theta_{i-2}-\theta_{i+1}}{\theta_{i-1}-\theta_i},
\qquad \qquad 
\frac{\theta^*_{i-2}-\theta^*_{i+1}}{\theta^*_{i-1}-\theta^*_i}
\end{eqnarray*}
are equal and independent of $i$ for $2 \leq i \leq D-1$
\cite[p.~263]{bannai}.
Denote this common value by $b+b^{-1}+1$ and assume that $b$ is
not a root of unity.
Further assume that, in the notation of
Bannai and Ito
\cite[p.~263]{bannai}, the given $Q$-polynomial structure
is type I
with $s\not=0$  and $s^*\not=0$.
Then for $b=q^2$ there exists a $\boxtimes_q$-action on the standard
module of $\Gamma$ for which 
the adjacency matrix acts as 
a $Z(T)$-linear combination of 
$1, x_{01}, x_{12}$ 
and the dual adjacency matrix acts as 
a $Z(T)$-linear combination of 
$1, x_{23}, x_{30}$.
\end{conjecture}

\begin{problem}
\rm A uniform poset \cite{uniform} is ranked and
has an algebraic structure
similar to that of a $Q$-polynomial distance-regular
graph. 
In \cite[p.~200]{uniform} 11 infinite families
of uniform posets are given.
For some uniform posets $P$
it might be possible to adapt the method of
the present paper to get 
an action of $\boxtimes_q$ on the standard module of $P$.
\end{problem}

\section{Acknowledgement} The authors would like to thank
Brian Curtin, Eric Egge, Mark MacLean, and Arlene Pascasio
for giving this manuscript a close reading
and offering many valuable suggestions.

\noindent Tatsuro Ito \hfil\break
\noindent Department of Computational Science \hfil\break
\noindent Faculty of Science \hfil\break
\noindent Kanazawa University \hfil\break
\noindent Kakuma-machi \hfil\break
\noindent Kanazawa 920-1192, Japan \hfil\break
\noindent email:  {\tt ito@kappa.s.kanazawa-u.ac.jp}

\bigskip

\noindent Paul Terwilliger \hfil\break
\noindent Department of Mathematics \hfil\break
\noindent University of Wisconsin \hfil\break
\noindent 480 Lincoln Drive \hfil\break
\noindent Madison, WI 53706-1388 USA \hfil\break
\noindent email: {\tt terwilli@math.wisc.edu }\hfil\break

\end{document}